\documentclass{amsart}

\newtheorem{thm*}[theorem]{Theorem}
\newtheorem{lem*}[theorem]{Lemma}
\newtheorem{prop*}[theorem]{Proposition}
\newtheorem{cor*}[theorem]{Corollary}
\newtheorem{conjecture*}{Conjecture}

\theoremstyle{definition}

\newtheorem{defn*}[theorem]{Definition}

\theoremstyle{remark}

\newtheorem{rem*}[theorem]{Remark}

\numberwithin{equation}{section}

\newcommand{\Tr}{\operatorname{Tr}}
\newcommand{\nc}{\newcommand}
\nc{\hM}{\widehat M}
\nc{\Da}{\Delta}
\nc{\da}{\delta}
\nc{\ta}{\theta}
\nc{\za}{\zeta}
\nc{\A}{\mathcal A}
\nc{\Om}{\Omega}
\nc{\Hj}{{\bar H}_{(2)}^j(\hM)}
\nc{\Hdot}{{\bar H}_{(2)}^{\bullet}(\hM)}
\nc{\Oj}{\Om_{(2)}^j(\hM)}
\nc{\vp}{\varphi}
\nc{\tildvp}{\tilde\vp}
\nc{\hK}{\widehat K}
\nc{\hN}{\widehat N}
\nc{\Cj}{C^j(\hK)}
\nc{\HjK}{{\bar H}_{(2)}^j(\hK)}
\nc{\HdotK}{{\bar H}_{(2)}^{\bullet}(\hK)}
\nc{\al}{\alpha}
\nc{\M}{{\mathbb H}}
\nc{\Tau}{T}

\def\al{\alpha}

\def\Om{\Omega}

\def\sup{\operatorname{sup}}

\def\Tr{\operatorname{Tr}}

\def\Det{\operatorname{Det}}

\def\A{{\mathcal A}}

\def\im{\operatorname{im}}

\def\<{\langle}
\def\>{\rangle}

\theoremstyle{remark}

\begin{document}


\title[Homotopy invariance of L$^2$ torsion for covering spaces]
{On the homotopy invariance of L$^2$ torsion for covering spaces.}
\author{Varghese Mathai}
\address{Department of Mathematics, University of Adelaide, Adelaide 5005,
Australia}
\email{vmathai@maths.adelaide.edu.au}
\author{Melvin Rothenberg}
\address{Department of Mathematics\\
University of Chicago, Chicago, U.S.A.}
\email{mel@math.uchicago.edu}
\thanks{Mel Rothenberg was supported in part by NSF Grant DMS 9423300.}

\date{MARCH 1997}


\begin{abstract}
We prove the homotopy invariance of $L^2$ torsion for covering spaces,
whenever the covering transformation group is either residually finite
or amenable.  In the case when the covering transformation group is
residually finite and when the $L^2$ cohomology of the covering space
vanishes, the homotopy invariance was established by L\"uck \cite{Lu}.
We also give some applications of our results.
\end{abstract}

\maketitle

\section*{Introduction}

$L^2$ analytic torsion for $L^2$ acyclic covering spaces with positive decay
was first studied in
\cite{M} and in \cite{L} and $L^2$ Reidemeister-Franz
torsion was first studied in \cite{CM}, see also \cite{LuR}, \cite{Lu}.
These $L^2$ torsion invariants were subsequently generalised by
using the theory of determinant lines of finitely
generated Hilbertian modules over finite von Neumann algebras,
which was initiated in \cite{CFM}.
Non-zero elements of the
determinant lines can be viewed as volume forms on the Hilbertian modules.
Using this, they constructed both $L^2$ combinatorial torsion and $L^2$
analytic torsion invariants
associated to flat Hilbertian bundles of determinant class
over compact polyhedra and manifolds,
as volume forms on the $L^2$ homology and $L^2$ cohomology
respectively, under the assumption that the covering space was of
determinant class.
These $L^2$ torsion invariants specialise to the Ray-Singer-Quillen
torsion and the classical Reidemeister-Franz
torsion respectively, in the finite dimensional case.
Using the results of  \cite{BFKM},
it was shown in \cite{CFM} that the combinatorial and analytic $L^2$ torsions
were equal
whenever the covering space is of determinant class. By this result
and calculations done in \cite{M} and \cite{L}, one obtains
a simplicial description of hyperbolic volume for a closed
3-dimensional hyperbolic manifold in terms of the $L^2$ Reidemeister-Franz
torsion.

In this paper, we prove the homotopy invariance of
$L^2$ torsion for covering spaces,
whenever the covering transformation group is either residually finite
or amenable.  In the case when the covering transformation group is
residually finite and when the $L^2$ cohomology of the covering space
vanishes, the homotopy invariance was established by L\"uck \cite{Lu}.
Using our main theorem, we can define the $L^2$ torsion
of a discrete group $\Gamma$, where $\Gamma$ is either residually finite
or amenable and
such that $B\Gamma$ is a finite CW complex. The $L^2$ torsion
of the discrete group $\Gamma$ is then an element of the
determinant line of the reduced $L^2$ cohomology of the
group $\Gamma$. We also give a new proof of  the classical result
that the hyperbolic volume of a closed 3-dimensional hyperbolic
manifold is a homotopy invariant, which actually only uses the part
of our theorem which was proved by L\"uck \cite{Lu}.
Our paper relies on the recent results that residually finite
covering spaces are of determinant class due to \cite{Lu},
\cite{BFK} and that amenable covering spaces are of determinant
class due to \cite{DM2}.

The authors thank Shmuel Weinberger for some helpful conversations.

\section{Preliminaries}

\subsection{Determinants and determinant lines}

Here we will review some of the results on the Fuglede-Kadison
determinant and on determinant lines of Hilbertian modules over finite
von Neumann algebras. The results of this section are mainly from
the paper \cite{CFM}. See also the papers \cite{CM}, \cite{FK}, \cite{Lu1}
and  \cite{LuR}.

Let $\A$ be a finite von Neumann algebra with a fixed faithful finite
normal trace $\tau :A\to\mathbb C$. Let $k$ denote the involution on $\A$.
Then $\A$ has a scalar product $\langle a,b\rangle=\tau (b*a)$ for
$a, b\in A$. Let $l^2(A)$ denote the completion of $\A$ with respect to
this scalar product.

A {\em Hilbert $\A$ module} is a Hilbert space  $M$ together with
a continuous left $\A$ module structure such that there is an
isometric $\A$-linear embedding of $M$ into $l^2(A)\otimes H$
for some Hilbert space $H$ ($\A$ acts trivially on $H$).
Note that the embedding is {\em not} part of the structure. $M$
is said to be {\em finitely generated} if $H$ can be chosen
above to be finite dimensional.

A Hilbert module has  a particular scalar product. If we ignore the scalar
 product, we get what we call a Hilbertian $\A$ module, more
precisely a {\em Hilbertian $\A$ module} is a topological
vector space $M$ with a left $\A$-module structure such that
there is a scalar product $\langle , \rangle$ on $M$ which generates
the topology of $M$, and such that $(M, \langle , \rangle )$ is a
Hilbert $\A$ module.

Any scalar product $\langle , \rangle $ on $M$ as above will be
called {\em admissible}. Suppose that $\langle , \rangle $ is
another scalar product on $M$. Then there is an operator
$A :M\to M$ such that $\langle v,w \rangle_1=\langle Av,w \rangle $
which satisfies
\begin{itemize}
\item[(1)]  $A$ is a linear homeomorphism since $\langle , \rangle$ and
$\langle , \rangle_1$ define the same topology,
\item[(2)]  self adjoint with respect to $\langle , \rangle$,
\item[(3)]  positive,
\item[(4)]  $A$ commutes with $\A$.
\end{itemize}
Any two admissible scalar products on a Hilbertian $\A$-module give rise to
isomorphic Hilbert $\A$-modules. In particular, if we choose an admissible
scalar product on $M$, it becomes a Hilbert $\A$ module and the {\em von
Neumann
dimension} $\dim_\tau (M)$, which is defined as the von Neumann
trace (defined in the next paragraph) of the orthogonal projection from
$\ell^2(\A)\otimes H$ onto $M$,
is independent of the choice of admissible scalar product.

The {\em commutant of $M$}, $B(M)$ is the algebra of all bounded
linear operators commuting with $\A$. Then there is a canonical
trace on $B(M)$, $\Tr_\tau$ defined as follows:
If $M$ were free, is isomorphic  to $l^2(A)\otimes{\mathbb C}^k$, then
$B(M)$ consists of $k\times k$ matrices with entries in $\A$,
acting on the right and for $f\in B(M)$,
$$
\Tr_{\tau}(f)=\sum_{i=1}^k \tau(f_{ii})
$$
is the trace. In the general case, $\Tr_{\tau}(f)= \Tr_{\tau}(i_M
\circ f\circ\pi_M)$, where $i_M$ denotes the embedding
of $M$ into a free module and $\pi_M$ denotes the
projection from the free module onto $M$.

Let $GL(M)$ denote the group of all {\em invertible
operators} in the commutant of $M$, $B(M)$. Then there
is a determinant function,
$$
\operatorname{Det}_\tau : GL(M)\to {\mathbb R}^+
$$
called the {\em Fuglede-Kadison determinant}, which satisfies
\begin{itemize}
\item[(1)]  $\operatorname{Det}_\tau(AB)=
\operatorname{Det}_\tau(A) \operatorname{Det}_\tau(B)$
\item[(2)]  $\operatorname{Det}_\tau(\lambda I)=|\lambda|^{\dim_{\tau} (M)}$
for $\lambda\in\mathbb C$
\item[(3)]  $\operatorname{Det}_{\lambda\tau}(A)=
\operatorname{Det}_\tau(A)^{\lambda}$ for $\lambda>0$
\item[(4)]  $\operatorname{Det}_\tau$ is continuous in the
norm topology on $GL(M)$
\end{itemize}
Its definition is as follows \cite{CFM}:
Let $A\in GL(M)$ and $A_t$ be a piecewise smooth path from
the identity to $A$ (since $GL(M)$ is connected, cf. \cite{Di}).
Define
\begin{equation}
\log\operatorname{Det}_\tau(A)=
\int_0^1\Re\Tr_\tau(A_t^{-1}\dot{A}_t) dt.
\tag{$*$}
\end{equation}
where $\Re$ denotes the real part. Using a result of Araki,
Smith and Smith \cite{ASS}, if $t\to A_t$ is a loop, then
$\int_0^1\Re \Tr_\tau(A_t^{-1}\dot{A}_t) dt=0$,
so that ($\ast$) is well defined.

Let $M$ be a finitely generated Hilbertian $\A$-module.
Define the {\em determinant line of} $M$, denoted by {$\det M$}, to be the
real vector space
generated by the symbols $\langle , \rangle$, one for
each {\em admissible} scalar product on $M$, with
the following relations : $\langle , \rangle_1$,
$\langle , \rangle_2$ are admissible scalar products on $M$, then
$$
\langle , \rangle_2=\operatorname{Det}_{\tau}(A)^{-1/2}
\langle , \rangle_1
$$
where $A\in GL(M)$, $A>0$ and
$$
\langle v, w\rangle_2=\langle Av, w\rangle_1
\quad\forall\ v, w\in M.
$$
It has been shown in \cite{CFM} that
$\det(M)$ is a real {\em one dimensional} vector space.
\nocite{*}

\subsection{$L^2$ torsion for covering spaces}
We will discuss $L^2$ torsion in the special case of normal covering
spaces; for the general case of Hilbertian $(\A-\pi)$ bimodules,
we refer to \cite{CFM}, \cite{BFKM}.

Let $\Gamma\to\hM\to M$ be a normal covering space of a compact
manifold. Let $K$ be a smooth triangulation of $M$, and $\hK$ be the induced
$\Gamma$-invariant triangulation of $\hM$. Then there is a natural
inner product on the space of finite support co-chains on
$\hK$, $C^j(\hK)$ where we declare the simplices to be an
orthonormal basis. Then $\Gamma$ acts by isometries on $\Cj$
and we consider the Hilbert space completion $C_{(2)}^j(\hK)$.
The coboundary operator
$$
d^K \, : \Cj\to C^{j+1}(\hK)
$$
induces a bounded operator
$$
d^{\widehat K} \, : C_{(2)}^j(\hK)\to C_{(2)}^{j+1}(\hK).
$$
Let $d^{\widehat K*}$ denote the $L^2$ adjoint of $d^{\widehat K}$, and the
combinatorial
Laplacian
$$
\Da_j^{\widehat K}=d^{\widehat K}d^{\widehat K*}+d^{\widehat K*}d^{\widehat
K}\, :
C_{(2)}^j(\hK)\to C_{(2)}^j(\hK)
$$
is a bounded self adjoint operator, which commutes with the
$\Gamma$-action, since $d^{\widehat K}$ and $d^{\widehat K *}$ do. Let
$\Da_j^{\widehat K+}$ denote
the restriction of $\Da_j^{\widehat K}$ onto the orthogonal complement
of $\operatorname{ker}\Da_j^{\widehat K}$. Assume now that
$\hM$ is of {\em determinant class}, that is, the Fuglede-Kadison determinant
of $\Da_j^{\widehat K+}$, $\operatorname{Det}_\tau(\Da_j^{\widehat K+})$ is
defined and is strictly positive for all $j\ge 0$, \cite{CM}, \cite{BFKM}.
It has been shown that
$\hM$ is of determinant class whenever $\Gamma$ is either residually finite
(cf. \cite{Lu}, \cite{BFK}) or amenable (cf. \cite{DM2}).

Since $\operatorname{ker}\Da_j^{\widehat K}\subset C_{(2)}^j(\hK)$, it
inherits an
admissible scalar product. By the combinatorial analogue of the
Hodge theorem, one sees that $\operatorname{ker}\Da_j^{\widehat K}$
is isomorphic to the reduced $L^2$ cohomology $\HjK$,
it follows that $\HjK$ gets an admissible
scalar product, i.e. an element in the determinant line $\det(\HjK)$. The
alternating
tensor product of these elements is an element in $\det(\HdotK)$,
which we denote by
$$
\hat\vp=\hat\vp(K)\in\det(\HdotK).
$$
The $L^2$ {\em Reidemeister-Franz torsion}, $\vp_{\hK}\in
\det(\HdotK)$ is defined as in \cite{CFM} to be
$$
\vp_{\hK}=\prod_{j=0}^{n}
\operatorname{Det}_\tau(\Da_j^{\widehat K +})^{\frac{(-1)^{j}j}{2}}\hat\vp(K).
$$
\begin{thm*}[Topological invariance of L$^2$ Reidemeister-Franz torsion
\cite{CFM}].
\\
Assume that $\hM$ be of determinant class. Then
$\vp_{\hK}$ is a topological invariant of $M$, whenever $\dim M$ is odd.
\end{thm*}
\begin{rem*}
When $\HjK=0 \;\;\forall j\geq0$, this result was first obtained by \cite{CM}.
In fact, $\vp_{\hK}$ is an invariant of the simple homotopy type of $K$, see
\cite{LuR} and 2.1 below. In the even dimensional case, the $L^2$
Reidemeister-Franz
torsion is not a topological invariant in general. Also, if $\dim M$ is
even and
$\HjK=0 \;\;\forall j\geq0$, then $\vp_{\hK}=1$.
\end{rem*}

There is also an analytic version of $L^2$ Riedemeister-Franz torsion
which can be defined in an analogous manner. It is
called $L^2$ {\em analytic torsion} and denoted by $\vp_{\hM}\in\det(\Hdot)$.
We refer to \cite{CFM} for details of its definition.

\begin{thm*}[Metric invariance of L$^2$ analytic torsion \cite{CFM}]
Assume that $\hM$ be of determinant class. Then
$\vp_{\hM}$ is independent of the choice of Riemannian metric
$g$ on $M$ whenever $\dim M$ is odd.
\end{thm*}
\begin{rem*}
When $\Hj=0\,\, \forall j\geq0$, this result was first obtained by \cite{M}
and \cite{L}. In \cite{CFM}, an even more general result is obtained.
In the even dimensional case, the $L^2$ analytic torsion does in general
depend on the choice of Riemannian metric. Also, if $\dim M$ is even and
$\Hj=0\,\, \forall j\geq0$, then $\vp_{\hM}=1$.
\end{rem*}

Under the de Rham isomorphism, which identifies $\HjK$ and $\Hj$ for
all $j$ (see \cite{Do}), the results of \cite{BFKM} are used to
prove the following theorem in \cite{CFM}.

\begin{thm*}[Equality of L$^2$ torsions \cite{CFM}, \cite{BFKM}]
Assume that $\hM$ be of determinant class. Then
via the identification of determinant lines induced by the de Rham isomorphism,
the L$^2$ analytic torsion and L$^2$ Reidemeister-Franz torsion
are equal, that is, $\vp_{\hM}=\vp_{\hK}$, whenever the dimension of $M$ is
odd.
\end{thm*}

\section{Homotopy invariance}

Here we establish the homotopy invariance of $L^2$ torsion for residually
finite
covering transformation groups
and for amenable covering transformation groups. As mentioned before, the
case when the
$L^2$ cohomology vanishes and when the covering transformation group is
residually finite, this result is due to \cite{Lu}.

We first formulate the problem of the homotopy invariance of $L^2$ torsion
as follows.
Let $f:M \to N$ be a homotopy equivalence of compact manifolds, $\widehat
f: \widehat M
\to \widehat N$ be the induced homotopy equivalence of the normal
$\Gamma$-covering
spaces $\widehat M$ and $\widehat N$. Then $\widehat f$ induces an
isomorphism on
the reduced $L^2$ cohomology, ${\widehat f}^* : {\bar
H}_{(2)}^\bullet(\widehat N)\to
 {\bar H}_{(2)}^\bullet(\widehat M)$. By 2.3 in \cite{CFM}, ${\widehat
f}^*$ induces
an isomorphism ${\widehat f}^*_*$ on determinant lines, ${\widehat f}^*_* :
\det {\bar H}_{(2)}^\bullet(\widehat N)\to \det{\bar
H}_{(2)}^\bullet(\widehat M)$.
Thus a homotopy
equivalence of manifolds induces a canonical isomorphism of determinant
lines of
$L^2$ cohomology.

\begin{conjecture*}[Homotopy invariance of $L^2$ torsion]
Let $f: M \to N$ be a homotopy equivalence of compact odd dimensional
manifolds. Suppose that
$\hM$ is of determinant class (or equivalently if $\hN$ is of determinant
class).
Via the identification of determinant lines of $L^2$ cohomology of normal
$\Gamma$
covering spaces as above, we conjecture that
$$
\phi_{\widehat M} = \phi_{\widehat N} \in  \det{\bar
H}_{(2)}^\bullet(\widehat M)
$$
\end{conjecture*}

This conjecture was first
formulated by L\"uck \cite{Lu} in the case when the $L^2$ cohomology
of $M$ is trivial.

Our first result will be a general formula, which generalises the ones in
\cite{LuR}, \cite{Lu} and \cite{Lu1}, in the case when the $L^2$ cohomology
is trivial. Let $A\in GL(n, {\mathbb Z}[\Gamma])$ be an invertible matrix
with entries in ${\mathbb Z}[\Gamma]$. Then it defines a bounded invertible
operator
$R_A : \ell^2(\Gamma)^n \to \ell^2(\Gamma)^n$ which commutes with the left
action
of $\Gamma$ on $\ell^2(\Gamma)^n$. So it has a positive Fuglede-Kadison
determinant
$\Det(R_A)$. Also $A$ defines an element in the Whitehead group of
$\Gamma$, and
it can be shown that the Fuglede-Kadison determinant $\Det$ induces a
homomorphism
$$
\Phi_\Gamma : Wh(\Gamma) \to {\mathbb R}^+
$$
which was defined in \cite{LuR} and \cite{Lu}. Now a homotopy equivalence
$f:M\to N$ defines an acyclic complex of cochains over the group ring of
$\Gamma$, $C^\bullet (M_f)$, where we now choose cell
decomposition for $M$ and $N$, and $f$ to be a cellular homotopy equivalence
and $M_f$ is the cellular mapping cone.
Then as in \cite{Mi}, this acyclic complex defines the Whitehead torsion
$\Tau(f)
\in Wh(\Gamma)$ of the homotopy equivalence $f$.
Note that via the identification above, one has
$\phi_{\widehat M} \otimes\phi_{\widehat N}^{-1} \in \mathbb R$.
Then we have the following

\begin{prop*}
 Let $f: M \to N$ be a homotopy equivalence of compact odd dimensional
manifolds. Suppose that $\hM$ is of determinant class.
Via the identification of determinant lines of $L^2$ cohomology of normal
$\Gamma$-covering spaces as above, one has
$$
\phi_{\widehat M}\otimes \phi_{\widehat N}^{-1} = \Phi_\Gamma(\Tau(f)) \in
{\mathbb R}^+
$$
\end{prop*}

\begin{rem*}
Farrell and Jones \cite{FJ} have proved that the
Whitehead group $Wh(\Gamma)$ is trivial
whenever $\Gamma$ is the fundamental group of a manifold with non-positive
sectional curvature. It follows from Proposition 2.1 that in this case, one has
$$
\phi_{\widehat M} = \phi_{\widehat N} \in  \det{\bar
H}_{(2)}^\bullet(\widehat M)
$$
i.e., the homotopy invariance of $L^2$ torsion,
provided the covering spaces are of determinant class.
Also Farrell and Jones \cite{FJ} conjectured that the
Whitehead group of any torsion-free
group vanishes, which would imply Conjecture 2 on the homotopy invariance
of $L^2$
torsion.

\end{rem*}

To prove this proposition, we recall some definitions and facts from \cite{CFM}
about abstract determinant class complexes. We begin with;

\begin{defn*} Consider a finitely generated Hilbertian module $M$ over $\A$.
A scalar product $\<\ ,\ \>$ on $M$ will be called {\it D-admissible}
if it can be represented in the form
$$\<v,w\>\ =\ \<A(v),w\>_1\quad\text{for}\quad v,w\in M,$$
where $\<\ ,\ \>_1$ is an admissible scalar product on $M$ and $A\in B(M)$
is an injective (possibly not invertible) homomorphism $A:M\to M$, which
is positive and self-adjoint with respect to $\<\ ,\ \>_1$, and
the following property is satisfied: if
$$A\ =\ \int_0^\infty \lambda dE_\lambda$$
is the spectral decomposition of $A$ and if
$$\phi(\lambda)\ =\ \dim_\tau(E_\lambda)\ =\ \Tr_\tau(E_\lambda)$$
denotes the corresponding spectral density function, then the integral
$$\int_0^\infty\ln(\lambda)d\phi(\lambda)\ >\ -\infty$$
is assumed to converge to a finite value.

 By a morphism of Hilbertian modules, we always mean one which
preserves the modules structure over the von Neumann algebra.
An injective homomorphism of Hilbertian modules $f: M\to N$ with dense image
will be called a {\em D-isomorphism} if for some,
and therefore all, admissible scalar product
$<\  ,\  >_N$ on $N$, the induced scalar product
$< \ , \ >_M$ on $M$ given by $<v,m>_M
= <f(v), f(w)>_N $, is {\em D-admissible}. Note
that any $D$-isomorphism of Hilbertian modules
$f:M\to N$, induces an isomorphism of determinant lines
$f_* : \det(M)\to \det(N)$, defined above, cf. \cite{CFM}.

A sequence of Hilbertian modules and morphisms
$0\to M' \stackrel{\alpha}\to M \stackrel{\beta}\to M'' \to 0$ is called
{\em D-exact}
is $\alpha$ is a
monomorphism, $\im \alpha = \ker \beta$ and the map induced by $\beta$, \
$M/\ker\beta \to M''$ is a $D$-isomorphism.

Let
$$
0\to C^0 \stackrel{d}\to C^1 \stackrel{d}\to\cdots \stackrel{d}\to C^N\to 0
$$
be a cochain complex of finite length formed by finitely
generated Hilbertian $\A$ modules and bounded linear maps
commuting with the action of $\A$. Here $\A$ is a finite von
Neumann algebra, such as the group von Neumann algebra of
$\Gamma$. Let $Z^i$ denote the submodule of cocycles and $B^i$
the submodule of coboundaries. Then the cochain complex $C^\bullet$
is said to be of {\em determinant class} if if the following
sequence
$$
0\to Z^i \stackrel{\alpha_i}\to C^i \stackrel{\beta_i}\to \bar{B^i}\to 0
$$
is $D$-exact, that is $im (\alpha_i) = ker (\beta_i)$
and the map induced by $\beta_i$, $C^i/ ker(\beta_i) \to \bar{B^i}$
is a $D$-isomorphism. This agrees with the notion of determinant class
discussed in section 1.2, (cf. \cite{CFM}).
\end{defn*}

We note as in \cite{CFM} that the property that a cochain complex
$C^\bullet$ of
finitely generated Hilbertian modules to be of determinant class depends only
on the homotopy type of $C^\bullet$ in the category of finitely generated
complexes of Hilbertian modules over $\A$.

\begin{defn*} Let $C^\bullet$ be a cochain complex which
is of {determinant class}. Then by Proposition 3.10 \cite{CFM}
there is a canonical isomorphism between the
determinant lines of the graded Hilbertian modules,
$$
\phi_C : \det(C) \to \det(\bar{H}(C))
$$
where $\bar{H}(C)$ denotes the reduced $L^2$ cohomology. This gives a
definition of
the $L^2$ {\em torsion} of an abstract cochain complex $C^\bullet$ which
is of determinant class,
$$
\phi_C \in  \det(C)^{-1} \otimes\det(\bar{H}(C))
$$
\end{defn*}

Let $K$ be a smooth triangulation on $M$, and $\hK$ be the induced
$\Gamma$-invariant triangulation on $\hM$.
It can be shown that
the $L^2$ torsion of the $L^2$ cochain
complex $C^\bullet_{(2)}(\widehat K)$ of $\widehat K$
is equal to the $L^2$  Reidemeister-Franz torsion of $\widehat K$,
(see Proposition 3.11 in \cite{CFM}).

\begin{prop*} Let $f: C \to C'$ be a cochain homotopy
equivalence of determinant class complexes. Then
$$
\phi_{C'}^{-1}\otimes \phi_C = \phi_{C_f}\in det(C_f)^{-1}.
$$
where $C_f$ denotes the mapping cone complex.
\end{prop*}

\begin{proof}
There is a short exact sequence
$$
0\to {C'}^{- 1} \to C_f \to C\to 0
$$
where $C_f$ denotes the mapping cone complex, and
${C'}^{-1}$ denotes the complex $C'$ shifted in degree
by $-1$. By Proposition 3.5 \cite{CFM}, one has the
canonical isomorphisms
$$
\det({C'}^{-1}) \otimes \det(C) \to \det(C_f)
$$
and
$$
\det(\bar{H}({C'}^{-1})) \otimes \det(\bar{H}(C)) \to \mathbb C
$$
since $C_f$ is acyclic. Therefore we get the canonical isomorphism
$$
\det(C')\otimes\det(\bar{H}(C'))^{-1} \otimes
\det(C)^{-1}\otimes\det(\bar{H}(C))
\to \det(C_f)^{-1}
$$
which yields the identity
$$
\phi_{C'}^{-1}\otimes \phi_C = \phi_{C_f}\in \det(C_f)^{-1}.
$$
\end{proof}

\begin{proof}[Proof of Proposition 2.1] By Proposition 2.5, it suffices to
determine the quantity $\phi_{C_f}\in \det(C_f)^{-1}$ in our particular case.
Firstly, since the complexes and the exact sequences
are based, we see that $\phi_{C_f}\in \mathbb R$.
Next, $C_f$ defines the Whitehead torsion $\Tau(f)\in Wh(\Gamma)$, and it is
clear that the $L^2$ torsion $\phi_{C_f}$ is just the Fuglede-Kadison
determinant of the Whitehead torsion $\Tau(f)$, $ \Phi_\Gamma(\Tau(f))$.
\end{proof}

Our next result is for amenable groups.

\begin{prop*}
Suppose that $\Gamma$ is a finitely presented amenable group. Then the
homomorphism
$$
\Phi_\Gamma : Wh(\Gamma) \to {\mathbb R}^+
$$
is trivial.
\end{prop*}

\begin{proof}
Let $a\in Wh(\Gamma)$. Then by \cite{Mi},
it is represented as the Whitehead torsion of a homotopy
equivalence $f:L\to K$ of finite CW complexes, which we
can assume without loss of generality
is an inclusion, i.e. $a = \Tau(f)$.
Let $\widehat K$ and $\widehat L$ denote the corresponding
$\Gamma$ normal covering complexes. The
relative cochain complex $C(\widehat K, \widehat L)$ is acyclic,
and so is its $L^2$ completion $C_{(2)}(\widehat K, \widehat L)$.
In particular, the combinatorial Laplacian
$\Delta_j^{\widehat K, \widehat L}$ is invertible, and
we see that
$$
\Phi_\Gamma(\Tau(f)) =
\prod_{j=0}^{n} \Det_\tau(\Delta_j^{\widehat K, \widehat L})^{\frac{(-1)^j
j}{2}} >0.
$$
Let ${\mathcal F}_K$ denote a fundamental domain for the action of the
group $\Gamma$ on $\widehat K$, and let ${\mathcal F}_L =
{\mathcal F}_K \cap {\widehat L}$, which is a fundamental domain for
the action of $\Gamma$ on $\widehat L$.

Since $\Gamma$ is amenable, by the F\o{}lner criterion for amenability
one gets a {\em regular exhaustion}, that
is, a sequence $\left\{X_m\right\}_{m=1}^{\infty}$ of {\em finite}
subsets of $\Gamma$ such that
$$
\lim_{m\to\infty}\frac{\sharp\partial_{\da} X_m}{\sharp X_m}=0
$$
where $\partial_{\da} X_m=\{\gamma\in\Gamma :
d(\gamma,X_m)<\da\;\;\text{and}\;\;
d(\gamma,\Gamma-X_m)<\da\}$ is a $\da$-neighbourhood of the "boundary" of
$X_m$, $d(\cdot,\cdot)$ denotes the word metric on $\Gamma$
and $ \sharp X_m$ denotes the number of elements in $X_m$.

One then gets {\em regular exhaustions}
$\big\{L_{m}\big\}^{\infty}_{m=1}$ of $\widehat L$,
and $\big\{K_{m}\big\}^{\infty}_{m=1}$
of $\widehat K$,
that is a sequences of finite subcomplexes of $\widehat L$ and
$\widehat K$ respectively such that

(1) $\displaystyle L_{m}=\bigcup_{g\in X_m} g.{\mathcal F}_L$ and
$\displaystyle K_m = \bigcup_{g\in X_m} g.{\mathcal F}_K$;

(2) $\displaystyle \widehat L=\bigcup^{\infty}_{m=1}L_{m}\;$
and $\displaystyle \widehat K=\bigcup^{\infty}_{m=1}K_{m}\;$

(3) Let $\dot{N}_{m,\delta}$ denote the number of elements $g\in \Gamma$
which have distance (with respect to the word metric in $\Gamma$) less than or
equal to $\delta$ from any element $g' \in \Gamma$ such that the
intersection of $g'.{\mathcal F}_L$
with  the topological boundary $\partial{L}_{m}$ of
$L_{m}$  is not empty, that is, $\dot{N}_{m,\delta} = \sharp\partial_{\da}
X_m$.
Let $N_m = \sharp X_m$, which also equals the number of translates
of  ${\mathcal F}_L$ which are contained in $L_m$. Then by the
amenability assumption, one has for every $\delta > 0$,
$$
\lim_{m\rightarrow\infty}\;\frac{{\dot{N}}_{m,\delta}}{N_{m}}=0.
$$
Note that $\dot{N}_{m,\delta}$ also equals the number of elements $g\in
\Gamma$
which have distance (with respect to the word metric in $\Gamma$) less than or
equal to $\delta$ from any element $g' \in \Gamma$ such that the
intersection of $g'.{\mathcal F}_K$
with  the topological boundary $\partial{K}_{m}$ of
$K_{m}$  is not empty, and $N_m$ also equals the number of translates of
${\mathcal F}_K$
which are contained in $K_m$.

We can then define the relative cochain complexes
$C(K_m, L_m)$ (note that $L_m = K_m \cap \widehat L$). Let
$\Delta_j^{(m)}$ denote the combinatorial Laplacian acting on
$C^j(K_m, L_m)$. Let $\Delta_j^{(m)+}$ denote the
combinatorial Laplacian acting on the orthogonal complement
of its null-space. Then by the results of section 4 in \cite{DM2},
which are easily generalised to the relative situation which is
being considered here, one has
$$
 \Det_\tau(\Delta_j^{\widehat K, \widehat L})\ge
\lim\sup_{m\to \infty}{\det}_{\mathbb C}(\Delta_j^{(m)+})^{\frac{1}{2N_m}} = 1
$$
since $\Delta_j^{(m)+}$ is an invertible matrix with integer entries.
To prove equality, we notice that if $E_m(\lambda)$ denotes
the spectral density function of  $\Delta_j^{(m)}$, then
\begin{equation}
\tag1 F_m(\lambda) = F_m(0) \;\;\;\;\mbox{for all}\;\;\; \lambda< K^{-2},
m\ge 0.
\end{equation}
where $F_m(\lambda) = \frac{E_m(\lambda)}{N_m}$ and $K$ is a constant as in
Lemma 2.2 in \cite{DM2} such that
$$
K^2\ge max\{||\Delta_j^{\widehat K, \widehat L}||, ||\Delta_j^{(m)+}||,
||{\Delta_j^{(m)+}}^{-1}||, ||{|\Delta_j^{\widehat K, \widehat L}}^{-1}||\}
$$
for all $m \ge 0$. Actually in Lemma 2.2 in \cite{DM2}, one obtains a
constant $K_1$ such that
$$
K_1^2\ge max\{||\Delta_j^{\widehat K, \widehat L}||, ||\Delta_j^{(m)+}||\}
$$
for all $m\ge 0$. Then an identical proof as in Lemma 2.2 in \cite{DM2}
gives the existence of a constant $K_2$ such that
$$
K_2^2\ge max\{||{\Delta_j^{(m)+}}^{-1}||,
||{|\Delta_j^{\widehat K, \widehat L}}^{-1}||\}.
$$
Finally, choose $K$ to be the greater of $K_1, K_2$.

The estimate $(1)$ says that $\Delta_j^{(m)}$ has a spectral gap at zero
which is bounded away from zero, independent of $m \ge 0$.
Let $\bar{F}(\lambda) = \lim\sup_{m\to\infty} F_m(\lambda)$, where we
assume without loss of generality that $\bar{F}(\lambda) $ is right continuous.
By Theorem 2.1 in \cite{DM2}, one has
$\bar{F}(\lambda) = F(\lambda)$, where
$F(\lambda)$ denotes the spectral density function of
$\Delta_j^{\widehat K, \widehat L}$. Therefore
$$
0 \le \bar{F}(\lambda) - \bar{F}(0) = \lim\sup_{m\to \infty} F_m(\lambda)
-\lim_{m\to\infty}F_m(0)
$$
$$
\le \sup\{|F_m(\lambda) - F_m(0)| : m\ge 0\} = 0
$$
if $\lambda < K^{-2}$, by (1).
Therefore one has,
\begin{equation}
\tag2 F(\lambda) =  \bar{F}(\lambda) =
\bar{F}(0)\;\;\;\;\mbox{for}\;\;\;\;\lambda<K^{-2}
\end{equation}
The estimate $(2)$ says that $\Delta_j^{\widehat K, \widehat L}$ also
has a spectral gap of size at least $K^{-2}$, at zero. Also observe that
$$
\int_0^{K^2} \sup\{F_m(\lambda) - F_m(0) : m\ge 0\} d\lambda
= \int_{K^{-2}}^{K^2} \sup\{F_m(\lambda) - F_m(0) : m\ge 0\} d\lambda <
\infty.
$$
By the dominated convergence theorem, one sees as in Lemma 3.3.1 in \cite{Lu}
that
$$
\int_0^{K^2} \frac{F(\lambda) - F(0)}{\lambda} d\lambda
= \lim_{m\to\infty}\int_{0}^{K^2} \frac{F_m(\lambda) - F_m(0)}{\lambda}
d\lambda.
$$
It follows that
$$
 \Det_\tau(\Delta_j^{\widehat K, \widehat L}) =
\lim_{m\to \infty}{\det}_{\mathbb C}(\Delta_j^{(m)+})^{\frac{1}{2N_m}} = 1.
$$
This completes the proof of the Proposition.
\end{proof}

Assembling Propositions 2.1 and 2.6 in this section, and using Theorem 0.5,
part a)
in \cite{Lu}, which states that
$$
\Phi_\Gamma : Wh(\Gamma) \to {\mathbb R}^+
$$
is trivial whenever $\Gamma$ is residually finite, and using Theorem 1.7 of the
previous section,
we obtain the following theorem
which gives significant evidence to the conjecture stated at the begining
of the
section.

\begin{thm*}
 Let $f: M \to N$ be a homotopy equivalence of compact, odd dimensional
manifolds. Via the identification of determinant lines of
$L^2$ cohomology of normal $\Gamma$
covering spaces as in the begining of the section,
and whenever $\Gamma$ is either residually finite
or $\Gamma$ is amenable, one has
$$
\phi_{\widehat M} = \phi_{\widehat N} \in  \det{\bar
H}_{(2)}^\bullet(\widehat M).
$$
\end{thm*}

\section{Applications}

In this section, we briefly give some applications of Theorem 2.7.

\subsection{$L^2$ torsion for discrete groups}
Our first application is the definition of the $L^2$ torsion for
a certain class of discrete groups.

\begin{prop*} Let $\Gamma$ be either a residually finite group or an
amenable group, whose classifying space $B\Gamma$
is a finite CW complex of odd dimension. Then one can
define the $L^2$ torsion of $\Gamma$ as
$$
\phi_\Gamma = \phi_{E\Gamma} \in \det({\overline H}^\bullet_{(2)}(\Gamma)).
$$
\end{prop*}

\begin{proof}
We note that since $\Gamma$ be either a residually finite group or an
amenable group such that $B\Gamma$ is a finite CW complex, it follows
that $E\Gamma$ is of determinant class, and therefore the $L^2$ torsion
$\phi_{E\Gamma}$ is well defined. Now $B\Gamma$ is only well defined upto
homotopy, and since both the determinant class condition depends only on
the homotopy type of $B\Gamma$,
and by Theorem 2.7, the $L^2$ torsion is a homotopy invariant, we therefore
see that the $L^2$ torsion $\phi_{E\Gamma}$ is independent of the choice
of $B\Gamma$, and depends only on the group $\Gamma$.
\end{proof}

We now conjecture the following.

\begin{conjecture*} Let $\Gamma$ be a discrete group whose classifying
space $B\Gamma$, is a finite CW complex of odd dimension.  Then one can
define the $L^2$ torsion of $\Gamma$ as
$$
\phi_\Gamma = \phi_{E\Gamma} \in \det({\overline H}^\bullet_{(2)}(\Gamma)).
$$
\end{conjecture*}

\subsection{Homotopy invariance of hyperbolic volume}
Our next application is a new proof of the homotopy invariance of the
hyperbolic volume of closed 3-dimensional hyperbolic manifolds. Although
it uses Theorem 2.7, it uses only the part that was proved by L\"uck \cite{Lu}.

\begin{prop*}
Let $M$ and $N$ be a homotopy equivalent closed 3-dimensional
hyperbolic manifolds. Then their hyperbolic volumes are equal.
\end{prop*}

\begin{proof} By the calculations in \cite{M} and \cite{L}, one sees that
the reduced $L^2$ cohomology of $M$ and $N$ vanish, and that the universal
covers
$\widetilde M =
\widetilde N = {\mathbb H}^3$, which is hyperbolic 3-space,
is of determinant class, and finally that the $L^2$ torsions of $\widetilde
M$ and
$\widetilde N$ can be calculated in terms of the hyperbolic volume,
$$- 6\pi\log\phi_{\widetilde M} = \mbox{vol}(M)$$
and
$$ - 6\pi\log\phi_{\widetilde N} = \mbox{vol}(N).$$
Now applying
Theorem 1.5 and Theorem 2.7, we see that
$$\mbox{vol}(M) = - 6\pi\log\phi_{\widetilde M} =
- 6\pi\log\phi_{\widetilde N} = \mbox{vol}(N).$$
\end{proof}


\end{document}